\documentclass[12pt,letterpaper]{amsart}
\usepackage{amssymb}

\newcommand{\indentalign}{\hspace{0.3in}&\hspace{-0.3in}}
\newcommand{\la}{\langle}
\newcommand{\ra}{\rangle}

\usepackage{amssymb}
\usepackage{amsthm}
\usepackage{amsxtra}
\newtheorem{theorem}{Theorem}
\newtheorem{remark}[theorem]{Remark}

\newtheorem{lemma}[theorem]{Lemma}

\numberwithin{equation}{section}
\numberwithin{theorem}{section}

\title[GWP for Zakharov and KGS]{Low regularity global well-posedness for the Zakharov and Klein-Gordon-Schr\"odinger systems}
\author{James Colliander}
\thanks{J.C. is partially  supported  by N.S.E.R.C. Grant RGPIN 250233-03 and the Sloan Foundation.}
\thanks{J.H. is supported by an NSF postdoctoral fellowship.}
\address{University of Toronto}
\author{Justin Holmer}
\address{University of California, Berkeley}
\author{ Nikolaos Tzirakis}
\address{University of Toronto}


\begin{document}

\begin{abstract}
We prove low-regularity global well-posedness for the 1d Zakharov system and 3d Klein-Gordon-Schr\"odinger system, which are systems in two variables $u:\mathbb{R}_x^d\times \mathbb{R}_t \to \mathbb{C}$ and $n:\mathbb{R}^d_x\times \mathbb{R}_t\to \mathbb{R}$.  The Zakharov system is known  to be locally well-posed in $(u,n)\in L^2\times H^{-1/2}$ and the Klein-Gordon-Schr\"odinger system is known to be locally well-posed in $(u,n)\in L^2\times L^2$.  Here, we show that the Zakharov and Klein-Gordon-Schr\"odinger systems are globally well-posed in these spaces, respectively, by using an available conservation law for the $L^2$ norm of $u$ and controlling the growth of $n$ via the estimates in the local theory. 
\end{abstract}

\maketitle

%\vfill
%\tableofcontents

%\newpage

\section{Introduction}
The initial-value problem for the one-dimensional Zakharov system is
\begin{equation}
\label{E:Zakharov}
\left\{
\begin{aligned}
&i\partial_tu + \partial_x^2 u = nu  \\
&\partial_t^2 n - \partial_x^2 n = \partial_x^2 |u|^2 \\
&u(x,0)=u_0(x), \; n(x,0)=n_0(x), \; \partial_tn(x,0)=n_1(x).
\end{aligned}
\right.
\end{equation}
Here $u:[0,T^*)\times \mathbb{R} \longmapsto \mathbb{R}, n:[0,T^*)
\times \mathbb{R} \longmapsto \mathbb{R}.$ This problem arises in
plasma physics.
Sufficiently regular solutions of \eqref{E:Zakharov} satisfy conservation of mass
\begin{equation}
\label{E:Zakmass}
M[u](t) = \int |u(t)|^2 \, dx=\int |u_{0}|^2 \, dx=M[u_0]
\end{equation}
and conservation of the Hamiltonian
\begin{align}
\label{E:Zakenergy}
H[u,n,\nu](t)& =\int \left( |\partial_x u(t)|^2 +\tfrac{1}{2}n(t)^2 +
  n(t)|u(t)|^2 +\tfrac{1}{2}\nu(t)^2 \right) \, dx \\
\notag &=H[u_{0},n_{0},\nu_{0}]
\end{align}
where $\nu(t)$ is such that $\partial_t n = \partial_x \nu$ and $\partial_t \nu = \partial_x (n+|u|^2)$.

The local-in-time theory in $X_{s,b}$ spaces has been established in
\cite{BC96}, \cite{GTV97}, the latter paper obtaining local
well-posedness (\textsc{LWP}) for the one-dimensional equation
\eqref{E:Zakharov} with  $(u_0,n_0,n_1)\in L^2\times H^{-1/2}\times
H^{-3/2}$ and for some more regular spaces $H^k\times H^s\times
H^{s-1}$ with various\footnote{The paper \cite{GTV97} actually gives a
  systematic treatment of \textsc{LWP} for higher dimensional versions
  of \eqref{E:Zakharov} as well.  Their result in dimension one uses
  the calculus techniques for obtaining $X_{s,b}$ bilinear estimates
  developed by Kenig, Ponce and Vega \cite{KPV96}, \cite{KPV96b}.  The only
  \textsc{LWP} result in \cite{GTV97} when $k=0$ is for dimension one,
  $s=-\frac{1}{2}$, and thus we have restricted exclusively to this
  case.} $k,s$.   As an immediate consequence of the local theory and \eqref{E:Zakmass},\eqref{E:Zakenergy}, one has global well-posedness (\textsc{GWP}) for $k=1$, $s=0$.  Pecher
\cite{Pec01}, using the low-high frequency decomposition method of
Bourgain \cite{Bou98}, proved \textsc{GWP} for $\frac{9}{10}<k<1$,
$s=0$.  This result was improved in \cite{Pec05} using the $I$-method of \cite{CKSTT02} to
obtain \textsc{GWP} for\footnote{Although the additional assumption
  $n_1\in \dot{H}^{-1}$ appears in the papers
  \cite{Pec01} and \cite{Pec05}, it can likely be removed by introducing
  suitable low frequency modifications to the energy identity; see
  \cite{GM94} Lemma A.1 p.\ 358.} $\frac{5}{6}<k<1$, $s=k-1$.  The
preceding GWP results are all based on the conservation of the
Hamiltonian \eqref{E:Zakenergy} or certain variants of the Hamiltonian. In this
paper, we prove \textsc{GWP} for $k=0$, $s=-\frac{1}{2}$, using a
scheme based on mass conservation \eqref{E:Zakmass} and subcritical slack in certain multilinear estimates at this regularity threshold.  In \cite{Hol05}, it is shown that the one-dimensional \textsc{LWP} theory of \cite{GTV97} is effectively sharp by adapting techniques of \cite{Bou93} and \cite{CCT03}. 
Thus, we establish \textsc{GWP} in the largest space for which \textsc{LWP} holds.
\begin{theorem}
\label{T:Zakglobal}
The Zakharov system \eqref{E:Zakharov} is globally well-posed for $(u_0,n_0,n_1)\in L^2\times H^{-1/2}\times H^{-3/2}$ and the solution $(u,n)$ satisfies \eqref{E:Zakmass} and
$$\|n(t)\|_{H_x^{-\frac{1}{2}}} + \|\partial_t n(t)\|_{H_x^{-\frac{3}{2}}} \leq \exp(c|t| \;\|u_0\|_{L^2}^2)\max(\|n_0\|_{H^{-\frac{1}{2}}}+\|n_1\|_{H^{-\frac{3}{2}}} \|u_0\|_{L^2}^2)
$$
\end{theorem}
%We remark that since this result is based only on \eqref{E:Zakmass} and estimates from the local theory, it applies equally well to \eqref{E:Zakharov} with $+nu$ replaced by $-nu$.  Interestingly, the corresponding Hamiltonian no longer (at least by itself) yields an \textit{a priori} bound on $\|u(t)\|_{H_x^1} + \|n(t)\|_{L_x^2}$ in this setting, and thus there is no known $H^1$ global result for nonlinearity $-nu$.

\begin{remark}
Since Theorem \ref{T:Zakglobal} is based on the mass conservation property \eqref{E:Zakmass} and the local theory, the same result applies to certain Hamiltonian generalizations of \eqref{E:Zakharov} for which global well-posedness was previously unknown. Indeed, if we write
$$
H [ u , n , \nu] (t) = \int | \partial_x u(t)|^2 + \frac{\alpha}{2} |n(t)|^2  + \frac{\beta}{2} |\nu (t) |^2 + \gamma  n(t) |u(t)|^2 dx
$$
and calculate
$$
\partial_t \left[ \begin{matrix} 
u \\
n \\
\nu 
\end{matrix}
\right]
= 
\left[
\begin{matrix}
-i & 0 & 0 \\
0 & 0 & \nabla \cdot \\
0 & \nabla & 0 
\end{matrix} 
\right]
\left[
\begin{matrix}
H_{\overline{u}} \\
H_n \\
H_\nu
\end{matrix}
\right]
$$
we find the evolution system
\begin{equation*}
\left\{
\begin{aligned}
& i\partial_t u + \partial_x^2 u = \gamma nu  \\
& \partial_t^2 n - \alpha \beta  \partial_x^2 n = \beta \gamma \partial_x^2|u|^2.
\end{aligned}
\right.
\end{equation*}
If we then choose $\alpha = \beta = -1$ and $\gamma = -1$, 
%\textbf{Justin: Assume you mean $\gamma$, not $\mu$?.  Also I think
%last line of above display should have $\beta \gamma |u|^2$ on
%right-hand side.  Please recheck.} 
we encounter a Hamiltonian evolution problem similar to \eqref{E:Zakharov} but with $+ nu$ replaced by $-nu$. The local theory for these problems coincides but the appearance of $\alpha = \beta = -1$ in the Hamiltonian $H$ precludes its use in obtaining a globalizing estimate.
\end{remark}

The initial-value problem for the $d$-dimensional Klein-Gordon-Schr\"odinger system with Yukawa coupling is
\begin{equation}
\label{E:KGS}
\left\{
\begin{aligned}
& i\partial_t u + \Delta u = -\gamma nu  && x\in \mathbb{R}^d, t\in \mathbb{R}\\
& \partial_t^2 n + \alpha \beta (1-\Delta)n = - \beta \gamma |u|^2  \\
& u(x,0)=u_0(x), \; n(x,0)=n_0(x), \; \partial_tn(x,0)=n_1(x).
\end{aligned}
\right.
\end{equation}
Here $\alpha, \beta, \gamma$ are real constants.
The solution satisfies conservation of mass
\begin{equation}
\label{E:KGSmass}
M[u](t)= \int |u(t)|^2 \,dx=\int |u_{0}|^2 \, dx=M[u_0]
\end{equation}
and conservation of the Hamiltonian
\begin{align}
\notag H[u,n,\nu](t) & = \int \left( |\nabla u(t)|^2 + \frac{1
  }{2\beta}|\partial_tn(t)|^2 + \frac{\alpha}{2}| \sqrt{-\Delta + 1} ~n(t)|
  ^2 + \gamma n(t) |u(t)|^2  \right) \, dx \\
\label{E:KGSenergy} & =H[u_{0},n_{0},\nu_{0}].
\end{align}

Pecher \cite{Pec04} proved that \eqref{E:KGS} is \textsc{LWP} for
$d=3$ in $L^2\times L^2\times H^{-1}$ and some more regular spaces
$H^k\times H^s\times H^{s-1}$ for various $k,s$ by following the
scheme developed for the Zakharov system in \cite{GTV97}.
 Provided that $\alpha >0$ and $\beta > 0$, energy conservation
 \eqref{E:KGSenergy} yields  \textsc{GWP} in the setting $1\leq d \leq
 3$ and $k=1$, $s=1$.  In the case when $\alpha > 0$ and $\beta >0$
 where the energy gives control on the $H^1$ norm, the low-high frequency separation method of
 Bourgain \cite{Bou98} has  been applied to \eqref{E:KGS} in
 \cite{Pec04} and the method of almost conservation laws of
 \cite{CKSTT02}   has been applied to \eqref{E:KGS} in \cite{Tzi05},
 to obtain \textsc{GWP} under the following assumptions: for $d=1$,
 $k=s$, $s>\frac{1}{2}$; for $d=2$, $k=s$, $s>\frac{\sqrt{17}-3}{2}$;
 for $d=3$, $k=s$, $s>\frac{7}{10}$; for $d=3$, $k,s>\frac{7}{10}$,
 $k+s>\frac{3}{2}$.  Moreover, in each of these cases, a polynomial in
 time bound is obtained for the growth of  the norms.  In this paper,
 we prove \textsc{GWP} for $d=3$, $k=s=0$, by a scheme involving
 \eqref{E:KGSmass} and direct application of the Strichartz estimates
 for the Schr\"odinger operator and Minkowski's integral inequality
 applied to the  Klein-Gordon Duhamel term.\footnote{Similar results
   hold for  $d=1$, $d=2$, although for expositional convenience, we
   have restricted  to the most delicate case $d=3$.}

\begin{theorem}
\label{T:KGSglobal}
The Klein-Gordon-Schr\"odinger system \eqref{E:KGS} in dimension $d=3$ is globally well-posed for $(u_0,n_0,n_1) \in L^2 \times L^2 \times H^{-1}$.  Moreover, the solution $(u,n)$ satisfies \eqref{E:KGSmass} and 
\begin{equation}
\label{E:504} \|n(t)\|_{L^2}+\|\partial_tn(t)\|_{H^{-1}} \leq \exp(c|t|\|u_0\|_{L^2}^2)\max((\|n_0\|_{L^2}+\|n_1\|_{H^{-1}}), \|u_0\|_{L^2}^2).
\end{equation}
\end{theorem}

\begin{remark}
In the case where $\alpha < 0$ and $\beta <0$, global well-posedness of \eqref{E:KGS} for
large smooth data was previously unknown. Since our proof of Theorem
\ref{T:KGSglobal} is based on the conservation of
$\| u(t) \|_{L^2}$, we do not require any Sobolev norm control obtained from the Hamiltonian and obtain global well-posedness for this case as well.
\end{remark}

The proof of both Theorem \ref{T:Zakglobal} and \ref{T:KGSglobal} apply essentially the same scheme, although invoke a different space-time norm in the local theory estimates.

\subsection{Outline of method}

We describe the globalization scheme for the Zakharov system and the
Klein-Gordon-Schr\"odinger system using the abstract
initial value problem posed at some time $t=T_j$
\begin{equation}
  \label{star}
  \left\{ 
\begin{aligned}
&Ku = F(u,n) \\
&Ln = G(u) \\
&(u,n)(T_j) = (u_j, n_j).
\end{aligned}
\right.
\end{equation}
Here $K$ and $L$ are linear differential operators of evolution type,
$F$ is a nonlinear term coupling the two equations together and $G$ is
a nonlinear term depending only upon $u$. The fact that $G$ does not
depend upon $n$ is used in our scheme.
Let $W(t) n_0$ denote the linear group $W(t)$ applied to initial data
$n_0$ solving the initial value problem $Ln = 0, n(0) =
n_0$. Similarly, let $S(t) u_0$ denote the solution of $Ku=0, u(0) =
u_0$. We denote with $W n_0 + L^{-1} g$ the solution of the linear initial
value problem $Lu = g, n(0) = n_0$. Similarly, $Su_0 + K^{-1} g$
denotes the solution of $Ku =g, u(0) = u_0$.

We solve the second equation in our system to define $n$ in terms of
the initial data $n_j$ and $u$ 
\begin{equation}
\label{nintermsofu}
n = W n_j + L^{-1} G(u)
\end{equation}
and insert the result into the solution formula for $u$ to obtain an
integrodifferential equation for $u$
\begin{equation}
\label{uduhamel}
u = Su_j + K^{-1} F( u, W n_j + L^{-1} G(u)).
\end{equation}
Local well-posedness for problems of the form \eqref{star} often
follows from a fixed point argument applied to \eqref{uduhamel}. The
fixed point analysis is carried out in a Banach space $X_{[T_j,
  T_{j-1}]}$ of functions defined on the spacetime slab $[T_j, T_{j+1}]
\times {\mathbb{R}}^d$.  The initial data are considered in function
spaces having the unitarity property with respect to the linear
solution maps
\begin{equation}
\label{unitary}
\| W(t) n_0 \|_{\mathcal{W}} = \| n_0 \|_{\mathcal{W}}, ~\| S(t) u_0
\|_{\mathcal{S}} = \| u_0 \|_{\mathcal{S}}, ~\forall~ t,
\end{equation}
and $X_{[T_j, T_{j+1}]} \subset C([T_j, T_{j+1}]; \mathcal{S}).$
For the applications we have in mind, the length of the time interval 
$\Delta_j := |[T_j, T_{j+1}]| $ is chosen to be small enough to prove a
contraction estimate and the smallness condition is of the form
\begin{equation*}
\Delta_j \leq \min( \| u_j \|_{\mathcal{S}}^{-\gamma},  \| n_j
\|_{\mathcal{W}}^{-\beta} )
\end{equation*}
for certain $\gamma, \beta > 0$.

Suppose that $\| u(t) \|_{\mathcal{S}} = \| u_0 \|_{\mathcal{S}}$ for
all times $t$ where solutions of \eqref{star} are well defined. If we
iterate the local well-posedness argument, we will have successive
time intervals $[T_j, T_{j+1}]$ with uniformly lower bounded lengths
unless $\| n_j \|_{\mathcal{W}}$ grows without bound as we increase
$j$. Suppose then at some time $T_j$ we have $\| n_j \|_{\mathcal{W}}
\gg \| u_j \|_{\mathcal{S}}^{\gamma/\beta}$ so that $\Delta_j = \| n_j \|_{\mathcal{W}}^{-\beta}.$

Since we have that \eqref{nintermsofu} and \eqref{unitary} hold, 
any growth in $\| n (t) \|_{\mathcal{W}}$
as $ t $ moves through the time interval $[T_j, T_{j+1}]$ is due to
the nonlinear influence of $u$ upon $n$ through the term $L^{-1} G(u)$.
Therefore, an estimate of the form
\begin{equation}
\label{incrementcontrol}
\| L^{-1} G(u) \|_{{L^\infty_{[T_j, T_{j+1}]} \mathcal{W}}} \leq
\Delta_j^\delta \widetilde{G} (\| u \|_{L^\infty_{[T_j, T_{j+1}]}
  \mathcal{S}} ) \ll \|n_j \|_{\mathcal{W}}
\end{equation}
permits an iteration of the local theory. Observe that the appearance of the conserved $\mathcal{S}$ norm of $u$
in this step suggests that we should retain this smallness property of
the $\mathcal{W}$ increment of $n$ over $[T_j, T_{j+1}]$ uniformly with
respect to $j$. 
%\textbf{[Justin:  I think you want $\|u\|_{L_{[T_j,T_{j+1}]}^\infty\mathcal{S}}
%$ inside the $\tilde G$ and then to mention use of the conservation
%law for $u$]}implies an almost conservation property for $n(t)$,
%\begin{equation}
%\label{abstractaclaw}
%\| n_{j+1} \|_{\mathcal{W}} \leq \| n_j \|_{\mathcal{W}} + \| L^{-1}
%G(u) (T_{j+1}) \|_{\mathcal{W}}. %< 2 \| n_j \|_{\mathcal{W}}.
%\end{equation}
%\textbf{[Justin: I think the $ < 2 \| n_j \|_{\mathcal{W}}$ will distract the reader into believing that we iterate with this estimate rather than the center term ]}
We then iterate the local well-posedness argument 
$$m = O \left( \frac{\| n_j \|_{\mathcal{W}}}{\Delta^\delta
    \widetilde{G} (\| u_0 \|_{\mathcal{S}})} \right)
$$
times with time steps of uniform size $\Delta = (2 \| n_j
\|_{\mathcal{W}})^{-\beta}$. 
%\textbf{[Justin: What is the point of
%starting at index $j$ and iterating forward by $m$ steps?  Shouldn't
%we start at $j=0$ and use $j$ as a dummy index over which to
%interate? Jim: I chose to write it this way since $j$ is the time
%when $\| n \|_{\mathcal{W}}$ is huge relative to $\| u
%\|_{\mathcal{S}}$. In principle either approach can be used. We can
%change it if you like. Perhaps Nikos can weigh in and we'll go with
%his choice.]} 
This extends the solution to the time
interval $[T_j, T_j + m \Delta]$ with 
$$
m\Delta = C(\| u_0 \|_{\mathcal{S}} ) \| n_j \|_{\mathcal{W}}^{1 - \beta + \delta \beta}.
$$
If $1 - \beta + \delta \beta \geq 0$, the scheme progresses to give
global well-posedness for \eqref{star}.

%\textbf{[Justin: I feel a shortcoming of the above outline is that it
%  does not advertise the point at which we appeal to the conservation
%  law for $u(t)$.  (It's buried in \eqref{incrementcontrol}) Jim: I
%  added a sentence highlighting the role of conservation near
%  \eqref{incrementcontrol}. I agree with Justin that this was not
%  sufficiently exposed.]}

Implementing this abstract scheme for specific systems requires a
quantification of the parameters $\beta$ and $\delta$ using the
local-in-time theory for the system. Notice that one way to force  $1 - \beta + \delta \beta \geq 0$ is by demanding $\beta \leq 1$.
 But $\beta$ is always bigger than $1$.  Still  
$1 - \beta + \delta \beta \geq 0$ can be greater or equal to zero
because of the contribution  of $\delta \beta$ term for certain
 $\beta>1$ and $\delta>0$. Calculations are required to obtain the
 parameters $\beta, \delta$ in any particular system. Unfortunately,
 we will often find that $\delta$ is very close to zero.  Thus the
 condition that $1 - \beta - \delta \beta \geq 0$
 fails to hold for many physical systems.  Nevertheless for the
 Klein-Gordon-Schr\"odinger  system the local well-posedness theory that we develop using the
 Strichartz's norms is sufficient for the above condition to hold. In
 particular we have  that $\beta=4$ and $\delta=\frac{3}{4}$ (see the
 proof of Theorem \ref{T:KGSglobal}) and thus $1 - \beta + \delta
 \beta=0$. This approach cannot be used for the Zakharov system.  The main reason is that
 the nonlinearity $G(u)$ has two derivatives (see equation
 \eqref{E:Zakharov}) and the local estimates are not as generous. The idea
 now is to perform the contraction argument for \eqref{uduhamel}  in a ball
$$B_{X_{[T_{j}^{\prime}, T_{j+1}^{\prime}]}}= \left\{
  u:\|u\|_{X_{[T_{j}^{\prime},  T_{j+1}^{\prime}]}} \leq
  (\Delta_{j}^{\prime})^{\alpha}\|u_{0}\|_{L^{2}} \right\} $$
where $\Delta_{j}^{\prime}=|[T_{j}^{\prime}, T_{j+1}^{\prime}]|$. This
idea is implemented here through the use
of the $X^{s,b}$ spaces with $b<1/2$. An easy consequence of this new
 iteration is  that the local time interval $\Delta_{j}^{\prime}$
 is larger or in other words (since $\Delta^{\prime}<1$) $\beta$ is smaller. In addition when we calculate the growth of the 
$\|n(t)\|_{\mathcal W}$ norm, it takes more time for this norm to
double in size.  In other words the new $\delta^{\prime}$ is bigger. Thus
since $m^{\prime}>m$ and $\Delta^{\prime}>\Delta$ we have a better chance to meet the requirement of 
$$m\Delta \gtrsim 1 \iff 1 - \beta + \delta \beta \geq 0.$$
The details are explained in the proof of Theorem \ref{T:Zakglobal}.
\begin{remark}
The abstract scheme described above can be applied to other evolution systems that have common features with systems \eqref{E:Zakharov} and
\eqref{E:KGS}. In particular, the scheme requires a satisfactory local
well-posedness theory in a Banach space that embeds in $C([T_j, T_{j+1}]; \mathcal{S})$, with $\| u(t) \|_{\mathcal{S}} = \| u_0 \|_{\mathcal{S}}$ holding true, and that the nonlinear term of the second equation is independent of $n$. As examples we mention the following systems.

The initial-value problem for the coupled Schr\"odinger-Airy equation
\begin{equation}
\label{E:SA}
\left\{
\begin{aligned}
&i\partial_t u + \partial_x^2u =\alpha un + \beta |u|^2u  && x\in \mathbb{R}, t\in \mathbb{R}\\
&\partial_tn + \partial_x^3n  = \gamma \partial_x|u|^2 \\
&u(x,0)=u_0(x), \; n(x,0)=n_0(x)
\end{aligned}
\right.
\end{equation}
This system arises in the theory of capillary-gravity waves.   The local well-posedness theory has been 
successively sharpened in \cite{BOP97}, \cite{BOP98}, \cite{CL05}, 
the last paper establishing 
local well-posedness for $(u_0,n_0)\in L^2\times H^s$ for $-\frac{3}{4}<s\leq -\frac{1}{2}$, and in some more regular spaces.
In \cite{Pec05b}, Pecher proved 
global well-posedness using $I$-method techniques for the harder Schr\"odinger-KdV system where the left hand side of the second equation of \eqref{E:SA} includes $n\partial_{x}n$, with $(u_0,n_0)\in H^s\times H^s$ and 
$s>\frac{3}{5}$ when $\beta=0$, and also for $s>\frac{2}{3}$ when $\beta\neq 0$ by dropping down from the $s=1$ 
setting in which conservation of energy yields global well-posedness.

Our scheme also applies to the Schr\"odinger-Benjamin-Ono system
\begin{equation}
\label{E:SBO}
\left\{
\begin{aligned}
&i\partial_t u + \partial_x^2u =\alpha un   && x\in \mathbb{R}, t\in \mathbb{R}\\
&\partial_tn + \nu \partial_x|\partial_x|n  = \beta \partial_x|u|^2 \\
&u(x,0)=u_0(x), \; n(x,0)=n_0(x)
\end{aligned}
\right.
\end{equation}
with $\alpha, \beta, \nu \in \Bbb R$. This system has been studied in
\cite{BOP98} where local well-posedness for 
$u_{0}\in H^{s}$ and $n_{0} \in H^{s-\frac{1}{2}}$, with $s \geq 0$
and $|\nu| \ne 1$ is established. In particular it is locally well-posed for
$(u_0,n_0)\in L^2\times H^{-1/2}$. Pecher proved \cite{Pec05c} global
well-posedness for $s>1/3$ under the parameter constraints $\nu >0$,\
$\frac{\alpha}{\beta}<0$ and also proved local well-posedness without
the restriction $|\nu| \ne 1$ but only for $s>0$. % We 
% believe that by following closely the above scheme and proving nonlinear estimates similar to these on Lemma \eqref{L:Zakmultest} we can 
% obtain global well-posedness for $(u_0,n_0)\in L^2\times H^{-1/2}$, for \eqref{E:SA} and for \eqref{E:SBO} with $|\nu| \ne 1$. Certain
%  extensions of these results we also believe are possible but we
%  explore them in a subsequent paper.    
In a forthcoming paper, we establish global well-posedness results for
\eqref{E:SA} and \eqref{E:SBO} with $|\nu| \neq 1$ for $(u_0 , n_0 )
\in L^2 \times H^{-1/2}.$
\end{remark}

\section{Basic estimates for the group and Duhamel terms}
Let $U(t)=e^{it\Delta}$ denote the free linear Schr\"odinger group.
For the 1d wave equation, it is convenient to factor the wave operator
$\partial_t^2-\partial_x^2 =
(\partial_t-\partial_x)(\partial_t+\partial_x)$, and work with
``reduced'' components, as was done in \cite{GTV97}.  Low frequencies
in the time-derivative initial data create some minor difficulties,
which we address in a manner slightly different than was done in \cite{GTV97}.  Consider an initial data pair $(n_0, n_1)$, and we look to solve $(\partial_t^2 - \partial_x^2) n =0$ such that $n(0)=n_0$, $\partial_t n(0)=n_1$.  Split $n_1 = n_{1L}+n_{1H}$ into low and high frequencies, and set $\hat{\nu}(\xi) = \frac{\hat{n}_{1H}(\xi)}{i\xi}$, so that $\partial_x \nu = n_{1H}$.  Let 
\begin{equation}
\label{E:Wdef}
\begin{aligned}
W_+(n_0,n_1)(t,x) &= \tfrac{1}{2}n_0(x-t) -\tfrac{1}{2}\nu(x-t) + \tfrac{1}{2} \int_{x-t}^x n_{1L}(y)\, dy\\
W_-(n_0,n_1)(t,x) &= \tfrac{1}{2}n_0(x+t) +\tfrac{1}{2}\nu(x+t) + \tfrac{1}{2} \int_x^{x+t} n_{1L}(y)\, dy
\end{aligned}
\end{equation}
so that 
$$
\begin{gathered}
(\partial_t\pm\partial_x)W_\pm(n_0,n_1)(t,x) = \tfrac{1}{2}n_{1L}(x) \\
W_\pm(n_0,n_1)(x,0)=\tfrac{1}{2}n_0(x)\mp \tfrac{1}{2}\nu(x) \\
\partial_tW_\pm(n_0,n_1)(x,0)=\mp\tfrac{1}{2}\partial_xn_0(x)+\tfrac{1}{2}\partial_x\nu(x)+\tfrac{1}{2}n_{1L}(x)
\end{gathered}
$$
and thus $n = W_+(n_0,n_1)+W_-(n_0,n_1)$ has the desired properties.  We shall also use the notation $W(n_0,n_1)=W_+(n_0,n_1)+ W_-(n_0,n_1)$.  Let 
$$G(t)(n_0,n_1) = \cos[t(I-\Delta)^{1/2}]n_0 + \frac{\sin[t(I-\Delta)^{1/2}]}{(I-\Delta)^{1/2}}n_1$$
be the free linear Klein-Gordon group, so that $(\partial_t^2+(1-\Delta))G(t)(n_0,n_1)=0$, $G(0)(n_0,n_1)=n_0$, $\partial_tG(0)(n_0,n_1)=n_1$.
Since our analysis involves tracking quantities whose size \textit{increments}, rather than \textit{doubles}, from one step to the next, it is imperative that we be precise about the definition of the following Sobolev norms.  When we write the norm $H^s$, we shall mean exactly
$$\| f\|_{H^s} = \left( \int_\xi (1+|\xi|^2)^s |\hat f(\xi)|^2 \, d\xi \right)^{1/2}.$$
Define the norm
$$\| f \|_{A^s} = \left( \int_{|\xi|\leq 1} |\hat f(\xi)|^2\, d\xi  + \int_{|\xi|\geq 1} |\xi|^{2s}|\hat f(\xi)|^2 \, d\xi \right)^{1/2}$$
Of course, $\|f\|_{A^s} \sim \|f\|_{H^s}$. Let
\begin{equation}
\label{E:Wnormdef}
\|(n_0,n_1)\|_{\mathcal{W}} = \left( \|n_0\|_{A_x^{-1/2}}^2 + \|n_1\|_{A_x^{-3/2}}^2 \right)^{1/2}.
\end{equation}
When working with a function of $t$, we use the shorthand $\|n(t)\|_\mathcal{W} = \|(n(t), \partial_t n(t))\|_{\mathcal{W}}$.  In our treatment of the Zakharov system, we shall track the size of the wave component $n(t)$ in the above norm.  Let
\begin{equation}
\label{E:Gnormdef}
\|(n_0,n_1)\|_{\mathcal{G}} = \left( \|n_0\|_{L_x^2}^2 + \| n_1\|_{H_x^{-1}}^2 \right)^{1/2}.
\end{equation}
Again, for functions of $t$, we use the shorthand $\|n(t)\|_\mathcal{G} = \|(n(t),\partial_tn(t))\|_{\mathcal{G}}$.  In our treatment of the Klein-Gordon-Schr\"odinger system, we shall track the size of the wave component $n(t)$ in the above norm.

In our treatment of the Zakharov system, we shall need to work in the Bourgain spaces.  We define the Schr\"odinger-Bourgain space $X^S_{0,\alpha}$, $\alpha\in \mathbb{R}$, by the norm
$$\|z\|_{X^S_{0,\alpha}} = \left( \iint_{\xi,\tau} \la \tau + |\xi|^2 \ra^{2\alpha} |\hat{z}(\xi,\tau)|^2 \, d\xi \; d\tau \right)^{1/2},$$
and the one-dimensional reduced-wave-Bourgain spaces $X^{W\pm}_{-\frac{1}{2},\alpha}$, for $\alpha\in \mathbb{R}$, as
$$\|z\|_{X^{W_\pm}_{-\frac{1}{2},\alpha}} = \left( \iint_{\xi,\tau} \la \xi \ra^{-1} \la \tau \pm \xi \ra^{2\alpha} |\hat{z}(\xi,\tau)|^2 \, d\xi \; d\tau \right)^{1/2}.$$

Let $\psi \in C_0^\infty (\mathbb{R})$ satisfy $\psi(t)=1$ on $[-1,1]$ and $\psi(t)=0$ outside of $[-2,2]$.  Let $\psi_T(t)= \psi(t/T)$, which will serve as a time cutoff for the Bourgain space estimates.  For clarity, we write $\psi_1(t) = \psi(t)$.  The following two lemmas are standard in the subject, although we are focusing attention particularly on the exponent of $T$ in these estimates.

\begin{lemma}[Group estimates] \label{L:Group}
Suppose $T\leq 1$.
\begin{enumerate}
\item \label{I:GroupSch}\emph{Schr\"odinger}.  $\|U(t)u_0 \|_{C(\mathbb{R}_t; L_x^2)} = \|u_0\|_{L_x^2}$.\\
 If $0\leq b_1\leq\tfrac{1}{2}$, then $\|\psi_T(t) U(t)u_0\|_{X_{0,b_1}^S} \lesssim T^{\frac{1}{2}-b_1}\|u_0\|_{L^2}$.  \\
\emph{(Strichartz Estimates)}.  If $2\leq q \leq \infty$, $2\leq r \leq \infty$, $\frac{2}{q}+\frac{d}{r}=\frac{d}{2}$, excluding the case $d=2$, $q=2$, $r=\infty$, then $\|U(t)u_0\|_{L_t^qL_x^r} \lesssim \|u_0\|_{L^2}$
\item \label{I:Groupwave}\emph{1-d Wave}.  $\|W(t)(n_0,n_1)\|_{C([0,T]; \mathcal{W}_x)} \leq (1+T) \|(n_0,n_1)\|_{\mathcal{W}}$.\\
If $0\leq b \leq \frac{1}{2}$, $\|\psi_T(t) W_\pm(t)(n_0,n_1)\|_{X^{W_\pm}_{-\frac{1}{2},b}} \lesssim T^{\frac{1}{2}-b}\|(n_0,n_1)\|_{\mathcal{W}}$.
\item \label{I:GroupKG}\emph{Klein-Gordon}.  $\| G(t)(n_0,n_1)\|_{C(\mathbb{R}_t; \mathcal{G}_x)} = \|(n_0,n_1)\|_{\mathcal{G}}$.
\end{enumerate}
\end{lemma}
\begin{remark}
It is important that the first estimate in \eqref{I:Groupwave} and the identity in \eqref{I:GroupKG} do not have implicit constant multiples on the right-hand side, as these estimates will be used to deduce almost conservation laws.  The $(1+T)$ prefactor in the first estimate of \eqref{I:Groupwave} arises from the low frequency terms.  Had we made the assumption that $n_1\in \dot{H}^{-1}$, this term could be removed and the norm $\mathcal{W}$ redefined so that equality is obtained.  The $(1+T)$ prefactor will not cause trouble in our iteration since $T$ will be selected so that $T\|(n_0,n_1)\|_\mathcal{W}$ functions as an increment whose size is on par with the increment arising from the Duhamel terms (see the proof of Theorem \ref{T:Zakglobal} for details).
\end{remark}
\begin{proof}
The Strichartz estimates quoted in \eqref{I:GroupSch} were established
in \cite{Str77} (for a more recent reference, see \cite{KT98}).   The first assertion in \eqref{I:GroupSch} is immediate by Plancherel's theorem.  For the second assertion in \eqref{I:GroupSch}, we note that $[\psi_T(t)U(t)u_0]\sphat(\xi,\tau) = (\psi_T)\sphat(\tau+\xi^2)\hat{u}_0(\xi)$, and consequently
$$\|\psi_T(t)U(t)u_0 \|_{X_{0,b_1}^S} \leq c \|\psi_T\|_{H^{b_1}}\|u_0\|_{L^2}$$To complete the proof of the estimate, we note that 
$$\|\psi_T\|_{H^{b_1}} \leq \|\psi_T\|_{L^2} + \|\psi_T\|_{\dot{H}^{b_1}} = T^{\frac{1}{2}}\|\psi_1\|_{L^2} + T^{\frac{1}{2}-b_1}\| \psi_1\|_{\dot{H}^{b_1}}$$
by scaling.    

For the first assertion in \eqref{I:Groupwave}, let $f(x,t)$ solve the linear wave equation 
\begin{equation}
\label{E:701}
\partial_t^2 f - \partial_x^2 f = 0
\end{equation}
with initial data $f(x,0)=n_0(x)$, $\partial_tf(x,0)=n_1(x)$.  Let $P_H$ be the projection onto frequencies $|\xi|\geq 1$, and $P_L$ be the projection onto frequencies $|\xi| \leq 1$.  Let $D^{-3/2}$ be the multiplier operator with symbol $|\xi|^{-3/2}$. By applying $D^{-3/2}P_H$ to \eqref{E:701}, multiplying by $D^{-3/2}P_H \partial_t f$ and integrating in $x$, we obtain the conservation identity
\begin{equation}
\label{E:703}
\|P_Hf(T)\|_{A_x^{-1/2}}^2 + \|\partial_t P_H f(T)\|_{A_x^{-3/2}}^2 = \|P_Hn_0\|_{A^{-1/2}}^2 + \|P_Hn_1\|_{A^{-3/2}}^2.
\end{equation}
To obtain low frequency estimates, we work directly from the explicit formula
\begin{equation}
\label{E:702}
f(x,t) = \tfrac{1}{2}n_0(x+t) + \tfrac{1}{2}n_0(x-t) + \tfrac{1}{2}\int_{x-t}^{x+t}n_1(y) \, dy.
\end{equation}
By applying $P_L$ and then directly estimating, we obtain
\begin{equation}
\label{E:704}
\|P_Lf(T)\|_{L_x^2} \leq \|P_Ln_0\|_{L^2} + T\|P_Ln_1\|_{L^2}.
\end{equation}
After applying $\partial_t$ to \eqref{E:702}, it can be rewritten as
$$\partial_t f(x,t) = \tfrac{1}{2}\int_{x-t}^{x+t}\partial^2n_0(y) \, dy + \tfrac{1}{2}n_1(x+t) + \tfrac{1}{2}n_1(x-t).$$
Applying $P_L$ and then directly estimating, we obtain
\begin{equation}
\label{E:705}
\|P_L \partial_t f(T)\|_{L_x^2} \leq T\|P_Ln_0\|_{L^2} + \|P_Ln_1\|_{L^2}
\end{equation}
Combining \eqref{E:703}, \eqref{E:704}, and \eqref{E:705}, we obtain the claim.

The second part of \eqref{I:Groupwave} is proved similarly to the second part of \eqref{I:GroupSch}.

For \eqref{I:GroupKG}, let $f(t,x)$ solve
\begin{equation}
\label{E:512}
\partial_t^2f - \Delta f + f = 0
\end{equation}
with initial data $(f(0),\partial_tf(0)) = (n_0,n_1)$.  Let $E$ be the multiplier operator with symbol $(1+|\xi|^2)^{-1/2}$.  Apply $E$ to \eqref{E:512}, then multiply by $\partial_tEf$, and finally integrate in $x$ to obtain the asserted conservation law.
\end{proof}

Let 
$$U\ast_R z(t,x) = \int_0^t U(t-t')z(t',x) \, dt'$$
denote the Duhamel operator corresponding to the Schr\"odinger operator, so that $(i\partial_t + \Delta)U\ast_R z(t,x) = iz(t,x)$, $U\ast_R z(0,x)=0$.  Let
$$W_\pm\ast_R z(t,x) = \tfrac{1}{2}\int_0^t z(t-s,x\mp s) ds$$
so that 
\begin{gather*}
(\partial_t \pm \partial_x)W_\pm\ast_R z(t,x) = \tfrac{1}{2}z(t,x)\\
W_\pm \ast_R z(0,x)= 0, \quad \partial_tW_\pm\ast_R z(0,x) = \tfrac{1}{2}z(0,x).
\end{gather*}
It follows that if we set $n=W_+\ast_R z - W_-\ast_R z$, then $(\partial_t^2 - \partial_x^2)n =\partial_x z$ and $n(0,x)=0$, $\partial_tn(0,x) = 0$, so we define
\begin{equation}
\label{E:706}
W \ast_R z = W_+ \ast_R z - W_-\ast_R z.
\end{equation}
For the Klein-Gordon equation, let
$$G\ast_R z(t,x) = \int_0^t \frac{\sin[(t-t')(I-\Delta)^{1/2}]}{(I-\Delta)^{1/2}}z(t',x) \, dt'$$
so that  $(\partial_t^2+(I-\Delta))G\ast_R z = z$, $G\ast_R z(0,x)=0$, $\partial_t G\ast_R z(0,x)=0$.

\begin{lemma}[Duhamel estimates]
\label{L:Duhamel}
Suppose $T\leq 1$.
\begin{enumerate}
\item \label{I:DuhamelSch}\emph{Schr\"odinger}.  If $0\leq c_1< \frac{1}{2}$, then $\| U\ast_R z \|_{C([0,T]; L^2_x)} \lesssim T^{\frac{1}{2}-c_1}\|z\|_{X_{0,-c_1}^S}$.\\
If $0\leq c_1<\frac{1}{2}$, $0\leq b_1$, $b_1+c_1\leq 1$, then $\| \psi_T U\ast_R z \|_{X_{0,b_1}^S} \lesssim T^{1-b_1-c_1}\|z\|_{X_{0,-c_1}^S}$. 
\emph{(Strichartz Estimates)}. If $2\leq q \leq \infty$, $2\leq r \leq
\infty$, $\frac{2}{q}+\frac{d}{r}=\frac{d}{2}$ excluding the case
$d=2$, $q=2$, $r=\infty$, and similarly for $\tilde{q}$, $\tilde{r}$,
then $\|U\ast_R z\|_{C([0,T];L_x^2)}+ \|U\ast_R z
\|_{L_{[0,T]}^qL_x^r} \lesssim
\|z\|_{L_{[0,T]}^{\tilde{q}'}L_x^{\tilde{r}'}}$, where $'$ indicates
the H\"older dual exponent ($\tfrac{1}{p} + \tfrac{1}{p'} = 1$).
\item \label{I:Duhamelwave}\emph{1-d Wave}.  If $0\leq c< \frac{1}{2}$, then $\| W \ast_R z \|_{C([0,T]; \mathcal{W}_x)} \lesssim T^{\frac{1}{2}-c}\|z\|_{X^{W\pm}_{-\frac{1}{2},-c}}$.\\
If $0\leq c<\frac{1}{2}$, $0\leq b$, $b+c\leq 1$, then $\| \psi_T W_\pm \ast_R z \|_{X_{-\frac{1}{2},b}^{W\pm}} \lesssim T^{1-b-c}\|z\|_{X_{-\frac{1}{2},-c}^{W\pm}}$.
\item \label{I:DuhamelKG} \emph{Klein-Gordon}.  $\|G\ast_R z\|_{C([0,T]; \mathcal{G}_x)} \lesssim \|z\|_{L_{[0,T]}^1 H_x^{-1}}$.
\end{enumerate}
\end{lemma}

\begin{proof}
The second assertion in each of \eqref{I:DuhamelSch} and \eqref{I:Duhamelwave} is \cite{GTV97} Lemma 2.1(ii).  For the Stricharz estimates quoted in \eqref{I:DuhamelSch}, see \cite{Str77}\cite{KT98}.   We next establish the first part of \eqref{I:DuhamelSch}.    We begin by establishing the bound
\begin{equation}
\label{E:D3}
\left\| \psi_T(t)U \ast_R z(x,t) \right\|_{L_t^\infty L_x^2} \leq
cT^{\frac{1}{2}-c_1} \|f\|_{X_{0,-c_1}^S}. 
\end{equation}
Let $f_\xi(t) = e^{it\xi^2}\hat z(\xi,t)$, where $\hat{\;}$ denotes the Fourier transform in the $x$-variable only.  We have
\begin{align}
\indentalign
\notag  \left\| \psi_T(t)U \ast_R z(x,t) \right\|_{L_t^\infty L_x^2} \\
\notag &= \left\| \psi_T(t)\int_0^t f_\xi(t') dt' \right\|_{L_t^\infty L_\xi^2}\\
\label{E:D4} &\leq  \left\| \psi_T(t)\int_0^t f_\xi(t') dt'
\right\|_{L_\xi^2 L_t^\infty}. 
\end{align}
Below we shall show that for a function $f(t)$ of the $t$-variable alone, we have the estimate
\begin{equation}
\label{E:D1}
\left\| \psi_T(t)\int_0^t f(t') dt' \right\|_{L_t^\infty} \leq cT^{\frac{1}{2}-c_1}\|f\|_{H^{-c_1}}.
\end{equation}
Assuming this, then it follows from \eqref{E:D4} that
$$\left\| \psi_T(t)U \ast_R z(x,t) \right\|_{L_t^\infty L_x^2} \leq cT^{\frac{1}{2}-c_1}\| \, \|f_\xi\|_{H_t^{-c_1}} \|_{L_\xi^2} = cT^{\frac{1}{2}-c_1} \|f\|_{X_{0,-c_1}^S}$$
completing the proof of \eqref{E:D3}.  Now we show \eqref{E:D1}.  Break $f(t)=f_+(t)+f_-(t)$ where $\hat{f}_-(\tau) = \chi_{|\tau|<\frac{1}{T}}\hat{f}(\tau)$ and $\hat{f}_+(\tau) = \chi_{|\tau|>\frac{1}{T}} \hat{f}(\tau)$.  Then for $f_-$, we have
$$\left\| \psi_T(t) \int_0^t f_-(s) \, ds \right\|_{L_t^\infty} \leq T^{1/2}\|f_-\|_{L^2} \leq T^{\frac{1}{2}-c_1} \|f_-\|_{H^{-c_1}}.$$
We compute
\begin{align*}
\indentalign \int_0^t f_+(s) \, ds = (\chi_{[-t,0]}\ast f_+)(0) \\
&= \frac{1}{2\pi} \int ( \chi_{[-t,0]}\ast f_+)\sphat(\sigma) \,d\sigma \\
&= \frac{1}{2\pi} \int_\sigma \frac{1-e^{-it\sigma}}{i\sigma} \hat{f}_+(\sigma) \, d\sigma
\end{align*}
and hence
$$\left[ \psi_T(t)\int_0^t f_+(s)\, ds\right]\sphat(\tau) =
\int_\sigma \frac{T\hat\psi(T(\tau-\sigma))-T\hat \psi(T\tau)}{\sigma}
\hat{f}_+(\sigma) \, d\sigma .$$
Thus,
\begin{align*}
\indentalign \left\| \psi_T(t) \int_0^t f_+(s)\, ds \right\|_{L_t^\infty} \\
&\leq \|T\hat\psi(T\tau)\|_{L_\tau^1} \left\| \frac{ \hat{f}_+(\tau)}{\tau} \right\|_{L_\tau^1} \\
&\leq \left( \int_{|\tau|>\frac{1}{T}} \frac{d\tau}{|\tau|^{2-2c_1}} \right)^{1/2}\|f_+\|_{H^{-c_1}} \\
&\leq T^{\frac{1}{2}-c_1}\|f\|_{H^{-c_1}}
\end{align*}
establishing \eqref{E:D1}.  It remains only to show continuity, i.e.\ that for a fixed $z\in X_{0,-c_1}^S$ and each $\epsilon>0$, there is $\delta=\delta(\epsilon, z)>0$ such that if $|t_2-t_1|<\delta$, then
$$\|\psi_T U\ast_R z(x,t_2)- \psi_T U\ast_R z(x,t_1)\|_{L_x^2} <\epsilon.$$
By an $\epsilon/3$ argument appealing to \eqref{E:D3}, it suffices to establish this statement for $z$ belonging to the dense class $\mathcal{S}(\mathbb{R}^2) \subset X_{0,-c_1}^S$.  However, if $z\in \mathcal{S}(\mathbb{R}^2)$, we have $\partial_t (U \ast_R z) = z +i\Delta (U\ast_R z)$ and the fundamental theorem of calculus and \eqref{E:D3} imply that
$$\| U \ast_R z(\cdot, t_2) - U \ast_R z(\cdot, t_1)\|_{L_x^2} \leq c(t_2-t_1)( \|z\|_{L_t^\infty L_x^2} +\|\Delta z\|_{X_{0,-c_1}^S}). $$
The proof of the first assertion of \eqref{I:Duhamelwave} proceeds in analogy to the above proof, first establishing the bound
$$\|W_\pm \ast_R z \|_{L_t^\infty H_x^{-1/2}} \lesssim T^{\frac{1}{2}-c}\|z\|_{X^{W\pm}_{-\frac{1}{2},-c}}.$$
The continuity statement is deduced by a density argument as in the previous paragraph, and finally the bound as stated on $W\ast_R z$ follows by the identity
$$\partial_t W\ast_R z = \partial_x W_+ \ast_R z + \partial_x W_-
\ast_R z .$$
The proof of \eqref{I:DuhamelKG} follows from an application of Minkowskii's integral inequality, with the continuity statement deduced by a density argument as in the previous paragraph.
\end{proof}

\section{1-d Zakharov system}
\label{S:1dZak}

In this section, we prove Theorem \ref{T:Zakglobal}.  We shall make use of the conservation law \eqref{E:Zakmass} to control the growth of $u(t)$ from one local time step to the next.  We track the growth of $n(t)$ in the norm $\mathcal{W}$ defined in \eqref{E:Wnormdef} using the estimates from the local theory.  We now state the needed estimates from the local theory of \cite{GTV97}.
 
\begin{lemma}[Multilinear estimates] \quad
\label{L:Zakmultest}
\begin{enumerate}
\item \label{I:ZakSchcoup}If $\tfrac{1}{4}<b_1,c_1,b<\frac{1}{2}$ and $b+b_1+c_1\geq 1$, then $\|n_\pm u \|_{X_{0,-c_1}^S} \lesssim \|n_\pm \|_{X_{-\frac{1}{2},b}^{W\pm}} \|u\|_{X_{0,b_1}^S}$.
\item \label{I:Zakwavecoup}If $\frac{1}{4}<b_1,c<\frac{1}{2}$ and $2b_1+c\geq 1$, then $\|\partial_x(u_1\bar{u}_2)\|_{X_{-\frac{1}{2},-c}^{W\pm}} \lesssim \|u_1\|_{X_{0,b_1}^S}\|u_2\|_{X_{0,b_1}^S}$.
\end{enumerate}
We remark that we can simultaneously achieve both optimal conditions $b+b_1+c_1=1$ and $2b_1+c=1$, for example by taking all four indices $b=b_1=c=c_1=\frac{1}{3}$.
\end{lemma}

\begin{proof}
\eqref{I:ZakSchcoup} is the case $k=0$, $\ell=-\frac{1}{2}$ in \cite{GTV97} Lemma 4.3, and \eqref{I:Zakwavecoup} is the case $k=0$, $l=-\frac{1}{2}$ in \cite{GTV97} Lemma 4.4.  The assumptions $b+c_1>\frac{3}{4}$, $b+b_1>\frac{3}{4}$ for Lemma 4.3 and $b_1+c>\frac{3}{4}$ for Lemma 4.4 that appear in \cite{GTV97} are not needed and we only have the requirements $b+b_1+c_1\geq 1$ for Lemma 4.3 and $2b_1+c\geq 1$ for Lemma 4.4.  The reason is that equation (4.30) in \cite{GTV97} on p.\ 424 is finite even if $\alpha_1<\frac{1}{2}$ since the range of integration is finite\footnote{This comment applies in the $k=0$, $\ell=-\frac{1}{2}$ setting but perhaps not in the general setting in which Lemmas 4.3,4.4 are stated.} (from $0$ to $\xi_1^2/4$).    Because relaxing this condition is essential to our method, we have included these proofs in the appendix so that they can be examined by the reader.
\end{proof}

\begin{proof}[Proof of Theorem \ref{T:Zakglobal}]  As discussed above, we can reduce the wave component $n=n_++n_-$ and recast \eqref{E:Zakharov} as
\begin{equation}
\label{E:Zakharovreduced}
\left\{
\begin{aligned}
&i\partial_tu + \partial_x^2 u = (n_++n_-)u  && x\in \mathbb{R}, t\in \mathbb{R}\\
&(\partial_t\pm \partial_x)n_\pm = \pm \tfrac{1}{2}\partial_x |u|^2 +\tfrac{1}{2}n_{1L} 
\end{aligned}
\right.
\end{equation}
which has the integral equation formulation
$$
\begin{aligned}
u(t) &= U(t)u_0 -iU\ast_R[(n_++n_-)u](t) \\
n_\pm(t) &= W_\pm(t)(n_0,n_1) \pm W_\pm\ast_R(\partial_x|u|^2)(t).
\end{aligned}
$$
Fix $0<T<1$, and consider the maps $\Lambda_S$, $\Lambda_{W\pm}$
\begin{align}
\label{E:301} \Lambda_S(u,n_\pm) &= \psi_TUu_0 + \psi_T U\ast_R[(n_++n_-) u] \\
\label{E:302} \Lambda_{W\pm}(u) &= \psi_TW_\pm(n_0,n_1) \pm \psi_T W_\pm \ast_R(\partial_x |u|^2).
\end{align}
We seek a fixed point $(u(t), n_\pm(t)) =
(\Lambda_S(u,n_\pm),\Lambda_{W\pm}(u))$.  
Estimating \eqref{E:301} in $X_{0,b_1}^S$, applying the first estimates in Lemma \ref{L:Group}\eqref{I:GroupSch}, \ref{L:Duhamel}\eqref{I:DuhamelSch} and following through with Lemma \ref{L:Zakmultest}\eqref{I:ZakSchcoup};  and estimating \eqref{E:302} in $X_{-\frac{1}{2},b}^{W_\pm}$, applying the first estimates in Lemma \ref{L:Group}\eqref{I:Groupwave}, \ref{L:Duhamel}\eqref{I:Duhamelwave} and following through with Lemma \ref{L:Zakmultest}\eqref{I:Zakwavecoup}, we obtain
\begin{align*}
\|\Lambda_S (u,n_\pm)\|_{X_{0,b_1}^S} & \lesssim T^{\frac{1}{2}-b_1}\|u_0\|_{L^2} + T^{1-b_1-c_1}\|n_\pm \|_{X_{-\frac{1}{2},b}^{W\pm}} \|u\|_{X_{0,b_1}^S} \\
\| \Lambda_{W\pm}(u)\|_{X_{-\frac{1}{2},b}^{W\pm}} & \lesssim T^{\frac{1}{2}-b}\|(n_0,n_1)\|_{\mathcal{W}} + T^{1-b-c} \|u\|_{X_{0,b_1}^S}^2
\end{align*}
and also
\begin{align*}
&\|\Lambda_S (u_1,n_{1\pm})-\Lambda_S (u_2,n_{2\pm})\|_{X_{0,b_1}^S} \\
& \qquad \lesssim  T^{1-b_1-c_1}(\|n_{1\pm} \|_{X_{-\frac{1}{2},b}^{W\pm}} \|u_1-u_2\|_{X_{0,b_1}^S}+\|n_{1\pm}-n_{2\pm} \|_{X_{-\frac{1}{2},b}^{W\pm}} \|u_2\|_{X_{0,b_1}^S}) \\
& \| \Lambda_{W\pm}(u_1)- \Lambda_{W\pm}(u_2)\|_{X_{-\frac{1}{2},b}^{W\pm}}  \lesssim T^{1-b-c} (\|u_1\|_{X_{0,b_1}^S}+\|u_2\|_{X_{0,b_1}^S})\|u_1-u_2\|_{X_{0,b_1}^S}.
\end{align*}
By taking $T$ such that 
\begin{gather}
\notag T^{\frac{3}{2}-2b_1-c_1}\|u_0\|_{L^2} \lesssim 1, \qquad T^{\frac{3}{2}-b-b_1-c}\|u_0\|_{L^2} \lesssim 1 \\
\label{E:306} T^{\frac{3}{2}-b-b_1-c_1}\|(n_0,n_1)\|_\mathcal{W} \lesssim 1 \\
\label{E:307} T^{\frac{3}{2}-2b_1-c}\|u_0\|_{L^2}^2 \lesssim \|(n_0,n_1)\|_\mathcal{W}
\end{gather}
%[{\bf Jim: It appears that we have incorporated an extra $T^{1/2 -
%  b_1}$ prefactor over what would be required to carry out the
%contraction argument. We should remark why we are doing this somewhere
%around here because the $3/2$ otherwise appears unnaturual.  Justin
%replies:  Do you think Nikos' addition to the introduction on this
%topic suffices and no further comment is necessary? Jim agress this
%is fine now.}]
one obtains sufficient conditions for a contraction argument yielding the existence of a fixed point $u\in X_{0,b_1}^S$, $n_\pm\in X_{-\frac{1}{2},b}^{W\pm}$ of \eqref{E:301}-\eqref{E:302} such that
\begin{equation}
\label{E:303}
\|u\|_{X_{0,b_1}^S} \lesssim T^{\frac{1}{2}-b_1}\|u_0\|_{L^2}, \qquad \|n_\pm \|_{X_{-\frac{1}{2},b}^{W\pm}} \lesssim T^{\frac{1}{2}-b}\|(n_0,n_1)\|_\mathcal{W}.
\end{equation}
Similarly estimating \eqref{E:301} in $C([0,T]; L_x^2)$ by applying Lemmas \ref{L:Group}\eqref{I:GroupSch},\ref{L:Duhamel}\eqref{I:DuhamelSch} and \eqref{E:303} shows that in fact $u\in C([0,T]; L_x^2)$. 
We may therefore invoke the conservation law \eqref{E:Zakmass} to conclude $\|u(T)\|_{L_x^2} = \|u_0\|_{L^2}$ and thus are concerned only with the possibility of growth in $\|n(t)\|_\mathcal{W}$ from one time step to the next.  Suppose that after some number of iterations we reach a time where $\|n(t)\|_\mathcal{W} \gg \|u(t)\|_{L_x^2}^2=\|u_0\|_{L_x^2}^2$.  Take this time position as the initial time $t=0$ so that $\|u_0\|_{L^2}^2 \ll \|(n_0,n_1)\|_{\mathcal{W}}$.  Then \eqref{E:307} is automatically satisfied and by \eqref{E:306}, we may select a time increment of size
\begin{equation}
\label{E:710}
T\sim \|(n_0,n_1)\|_{\mathcal{W}}^{-1/(\frac{3}{2}-b-b_1-c_1)} = \|(n_0,n_1)\|_{\mathcal{W}}^{-2}
\end{equation}
where the right-hand side follows by selecting the optimal condition $b+b_1+c_1=1$ in Lemma \ref{L:Zakmultest}\eqref{I:ZakSchcoup}.  Since
$$n= W(n_0,n_1) + W\ast_R (\partial_x |u|^2)$$
we can apply Lemma \ref{L:Group}\eqref{I:Groupwave}, \ref{L:Duhamel}\eqref{I:Duhamelwave} and follow through with \eqref{E:303} to obtain
\begin{align*}
\|n(T)\|_\mathcal{W} &\leq (1+T)\|(n_0,n_1)\|_\mathcal{W} + CT^{\frac{3}{2}-(2b_1+c)}\|u_0\|_{L^2}^2 \\
&\leq \|(n_0,n_1)\|_\mathcal{W} + CT^\frac{1}{2}(\|u_0\|_{L^2}^2+1)
\end{align*}
where $C$ is some fixed constant.  The second line above follows by selecting the optimal condition $2b_1+c=1$ in Lemma \ref{L:Zakmultest}\eqref{I:Zakwavecoup}, and using \eqref{E:710} to obtain $T\|(n_0,n_1)\|_{\mathcal{W}} \leq CT^\frac{1}{2}$.   From this we see that we can carry out $m$ iterations on time intervals, each of length \eqref{E:710}, where
\begin{equation}
\label{E:711}
m \sim \frac{\|(n_0,n_1)\|_\mathcal{W}}{T^\frac{1}{2}(\|u_0\|_{L^2}^2+1)}
\end{equation}
before the quantity $\|n(t)\|_\mathcal{W}$ doubles.  The total time we advance after these $m$ iterations, by \eqref{E:710} and \eqref{E:711}, is
$$mT \sim \frac{1}{\|u_0\|_{L^2}^2+1}$$
which is independent of $\| n(t) \|_{\mathcal{W}}$.

We can now repeat this entire procedure, each time advancing a time of
length $\sim (\|u_0\|_{L^2}^2+1)^{-1}$ (independent of the size of
$\|n(t)\|_\mathcal{W}$).  Upon each repetition, the size of
$\|n(t)\|_\mathcal{W}$ will at most double, giving the
exponential-in-time upper bound stated in Theorem \ref{T:Zakglobal}.

\end{proof}

\section{3-d Klein-Gordon Schr\"odinger system}

The goal of this section is to prove Theorem \ref{T:KGSglobal}.  For the Klein-Gordon-Schr\"odinger system \eqref{E:KGS}, no special multilinear estimates are needed. Instead, we will work in standard space-time norms and use Sobolev imbedding and the H\"older inequality.

We shall use the conservation law \eqref{E:KGSmass} to control the growth of $u(t)$ from one time step to the next and track the growth of $\|n(t)\|_\mathcal{G}$, where the $\mathcal{G}$ norm was defined in \eqref{E:Gnormdef}, by direct estimation.  For expositional convenience, we restrict to dimension $d=3$ although similar results do hold for $d=1$ and $d=2$.

\begin{proof}[Proof of Theorem \ref{T:KGSglobal}]
\eqref{E:KGS} has the integral equation formulation
\begin{align*}
 &u(t) = U(t)u_0 + iU\ast_R[nu](t) \\
 &n(t) = G(t)(n_0,n_1) + G\ast_R(|u|^2)(t).
\end{align*}
Define the maps $\Lambda_S$, $\Lambda_G$ as
\begin{align}
\label{E:501} &\Lambda_S(u,n) = Uu_0 + iU\ast_R[ nu] \\
\label{E:502} &\Lambda_G(u) = G(n_0,n_1) + G\ast_R(|u|^2).
\end{align}
Let $\text{Str}= L_{[0,T]}^{10/3}L_x^{10/3}\cap L_{[0,T]}^8L_x^{12/5}$.   We seek a fixed point $(u(t),n(t))=(\Lambda_S(u,n),\Lambda_G(u))$ in the space $[C([0,T];L_x^2)\cap \text{Str}]\times C([0,T]; \mathcal{G}_x)$.  Apply Lemma
\ref{L:Duhamel}\eqref{I:GroupSch} with $(\tilde{q},\tilde{r})=(\frac{20}{9},5)$ for $d=3$ to obtain
$$ \| \Lambda_S(u,n) \|_{C([0,T]; L_x^2)\cap \text{Str}} \lesssim \|u_0\|_{L^2} + T^{1/4}\|n\|_{L_{[0,T]}^\infty L_x^2} \|u\|_{L_{[0,T]}^\frac{10}{3}L_x^\frac{10}{3}}.$$
Estimate \eqref{E:502} in $C([0,T];\mathcal{G}_x)$ and apply Lemma \ref{L:Duhamel}\eqref{I:DuhamelKG} followed by Sobolev imbedding to obtain
\begin{equation}
\label{E:503}
\| \Lambda_G(u) \|_{C([0,T]; \mathcal{G}_x)} \leq \|(n_0,n_1)\|_\mathcal{G} +c
T^{3/4}\|u\|_{L_{[0,T]}^8L_x^{12/5}}^2
\end{equation}
where we estimated as: $\| |u|^2 \|_{L_{[0,T]}^1H_x^{-1}} \lesssim \| |u|^2 \|_{L_{[0,T]}^1L_x^{6/5}} \leq T^{3/4}\|u\|_{L_{[0,T]}^8L_x^{12/5}}^2$.
There are similar estimates for the differences $\Lambda_S(u_1,n_1)-\Lambda_S(u_2,n_2)$ and $\Lambda_G(u_1)-\Lambda_G(u_2)$.   If $T$ is such that
\begin{gather}
\notag T^{1/4}\|u_0\|_{L^2} \lesssim 1 \\
\label{E:511} T^{1/4}\|(n_0,n_1)\|_{\mathcal{G}}\lesssim 1 \\
\label{E:510} T^{3/4}\|u_0\|_{L^2}^2 \lesssim \|(n_0,n_1)\|_{\mathcal{G}}
\end{gather}
then a contraction argument implies there is a solution $(u,n)$ to \eqref{E:KGS} on $[0,T]$ such that
\begin{align}
\label{E:514} &\|u\|_{C([0,T];L_x^2)\cap \text{Str}} \lesssim \|u_0\|_{L^2} \\
\label{E:517} &\|n\|_{C([0,T]; \mathcal{G}_x)} \leq \|(n_0,n_1)\|_{\mathcal{G}} + cT^{3/4}\|u_0\|_{L^2}^2.
\end{align}
By the conservation of mass, we have $\|u(t)\|_{L_x^2} = \|u_0\|_{L^2}$, and are thus concerned only with the possibility that $\|n(t)\|_{\mathcal{G}}$ grows excessively from one local increment to the next.  Suppose that after some number of iterations $\|n(t)\|_{\mathcal{G}} \gg \|u(t)\|_{L^2}^2=\|u_0\|_{L^2}^2$.  Consider this time as the initial time so that $\|(n_0,n_1)\|_{\mathcal{G}} \gg \|u_0\|_{L^2}^2.$
Then \eqref{E:510} is automatically satisfied, and by \eqref{E:511}, we may thus take
\begin{equation}
\label{E:515}
T\sim \|(n_0,n_1)\|_\mathcal{G}^{-4}.
\end{equation}
We see from \eqref{E:517} that, after $m$ iterations, each of size \eqref{E:515}, where 
$$m\sim \frac{\|(n_0,n_1)\|_{\mathcal{G}}}{T^{3/4}\|u_0\|_{L^2}^2}$$
the quantity $\|n(t)\|_{\mathcal{G}}$ at most doubles. The total time advanced after these $m$ iterations is
$$mT \sim \frac{1}{\|u_0\|_{L^2}^2}.$$
We can now repeat this entire procedure, each time advancing a time of
length $\sim \|u_0\|_{L^2}^{-2}$ (independent of the size of
$\|n(t)\|_{\mathcal{G}}$).  Upon each repetition, the size of
$\|n(t)\|_{\mathcal{G}}$ will at most double, giving the
exponential-in-time upper bound stated in Theorem \ref{T:KGSglobal}. 

\end{proof}

\appendix
\section{Proof of the multilinear estimates (expository)}

In this section, we prove Lemma \ref{L:Zakmultest}.  The material here is taken from \cite{GTV97} Lemma 4.3, 4.4 with only a slight modification at one stage.  This modification was described in a note under the heading ``proof'' following the statement of Lemma \ref{L:Zakmultest}.  Given its importance in our scheme, the full proof is included here in detail.

We need the calculus lemmas:

\begin{lemma}[{\cite{GTV97}} Lemma 4.1]
\label{L:GTV4.1}
Let $f\in L^q(\mathbb{R})$, $g\in L^{q'}(\mathbb{R})$ for $1\leq q,q'\leq \infty$ and $\frac{1}{q}+\frac{1}{q'}=1$.  Assume that $f,g$ are nonnegative, even, and nonincreasing for positive argument.  Then $f\ast g$ enjoys the same properties.
\end{lemma}

Define
$$[\lambda]_+ = 
\left\{
\begin{aligned}
&\lambda  && \text{if }\lambda>0 \\
&\epsilon && \text{if }\lambda=0 \\
&0 && \text{if }\lambda<0 
\end{aligned}
\right.
$$

\begin{lemma}[{\cite{GTV97}} Lemma 4.2]
\label{L:GTV4.2}
Let $0\leq a_-\leq a_+$ and $a_++a_->\frac{1}{2}$.  Then $\forall \; s\in \mathbb{R}$,
$$\int_y \la y-s \ra^{-2a_+} \la y+s \ra^{-2a_-} \, dy \leq c\la s \ra^{-\alpha}$$
where $\alpha=2a_--[1-2a_+]_+$.
\end{lemma}

\begin{proof}[Proof of Lemma \ref{L:Zakmultest}\eqref{I:ZakSchcoup}]
We shall only do the $+$ case.   The estimate is equivalent to
$$|S| \leq c\|v\|_2\|v_1\|_2\|v_2\|_2$$
where
\begin{equation}
\label{E:600}
S = \int_\ast \frac{\hat{v} \hat{v}_1 \hat{v}_2 \la \xi \ra^{1/2}}{ \la \sigma \ra^b \la \sigma_1 \ra^{c_1} \la \sigma_2 \ra^{b_1}}
\end{equation}
with $\hat{v}=\hat{v}(\xi,\tau)$,  $\hat{v}_1=\hat{v}_1(\xi_1,\tau_1)$, $\hat{v}_2=\hat{v}_2(\xi_2,\tau_2)$, $\sigma_1 = \tau_1 + \xi_1^2$, $\sigma_2=\tau_2+\xi_2^2$, $\sigma = \tau+\xi$, and $*$ indicates the restriction $\xi_1=\xi+\xi_2$, $\tau_1=\tau+\tau_2$.  Indeed, for $\hat{v}_1 \in L^2$,
\begin{align*}
\indentalign \int_{\xi_1,\tau_1} \widehat{n_+ u}(\xi_1,\tau_1)  \la \sigma_1 \ra^{-c_1} \hat{v}_1(\xi_1,\tau_1) \, d\xi_1 d\tau_1 \\
&= \int_{\xi_1,\tau_1} \left[ \int_{\substack{\xi_1=\xi+\xi_2 \\
      \tau_1=\tau+\tau_2}} \hat{n}_+(\xi,\tau)
  \hat{u}(\xi_2,\tau_2)\right]  \la \sigma_1 \ra^{-c_1}
\hat{v}_1(\xi_1,\tau_1) d\xi_1 d\tau_1. 
\end{align*}
Let $\hat{v}(\xi,\tau) = \hat{n}_+(\xi,\tau) \la \xi \ra^{-1/2} \la \sigma \ra^b$ and $\hat{v}_2(\xi_2,\tau_2) = \hat{u}(\xi_2,\tau_2)\la \sigma_2 \ra^{b_1}$ to obtain \eqref{E:600}.

We note here that
\begin{align*}
\indentalign \sigma_1 -\sigma -\sigma_2 = \tau_1 + \xi_1^2 - (\tau+\xi) - (\tau_2+\xi_2^2) \\
&= \xi_1^2 - \xi-\xi_2^2 \\
&= (\xi_1-\tfrac{1}{2})^2 - (\xi_2-\tfrac{1}{2})^2
\end{align*}

In the case when $|\xi|\leq 1$, it suffices to estimate
\begin{align*}
\indentalign \int_{\xi_1,\tau_1} \int_{\substack{\xi_1=\xi+\xi_2 \\ \tau_1=\tau+\tau_2}} \frac{ \hat{v} \hat{v}_1 \hat{v}_2}{\la \sigma \ra^{1/4} \la \sigma_1 \ra^{1/4+} \la \sigma_2 \ra^{1/4+}} \\
&=\int_{\xi_1,\tau_1} \left( \frac{\hat{v}_1}{\la \sigma_1 \ra^{1/4+}} \right) \left( \int_{\substack{\xi_1=\xi+\xi_2 \\ \tau_1=\tau+\tau_2}} \frac{\hat{v} \hat{v}_2}{\la \sigma \ra^{1/4} \la \sigma_2 \ra^{1/4+}} \, d\xi_2 d\tau_2 \right) \, d\xi_1 d\tau_1 \\
&= \int_{x_1,t_1} \left[ \frac{\hat{v}_1}{\la \sigma_1 \ra^{1/4+}} \right]\spcheck \left[ \int_{\substack{\xi=\xi_1+\xi_2 \\ \tau=\tau_1+\tau_2}} \frac{\hat{v} \hat{v}_2}{\la \sigma \ra^{1/4} \la \sigma_2 \ra^{1/4+}} \, d\xi_2 d\tau_2 \right]\spcheck \\
&= \int_{x,t} \left[ \frac{\hat{v}_1}{\la \sigma_1 \ra^{1/4+}} \right]\spcheck \left[ \frac{\hat{v}}{\la \sigma \ra^{1/4}} \right]\spcheck  \left[ \frac{\hat{v}_2}{\la \sigma_2 \ra^{1/4+}} \right]\spcheck \\
&\leq \left\| \left[ \frac{\hat{v}_1}{\la \sigma_1 \ra^{1/4+}} \right]\spcheck \right\|_{L_t^{8/3}L_x^4} \left\| \left[ \frac{\hat{v}}{\la \sigma \ra^{1/4}} \right]\spcheck \right\|_{L_t^4L_x^2} \left\| \left[ \frac{\hat{v}_2}{\la \sigma_2 \ra^{1/4+}} \right]\spcheck \right\|_{L_t^{8/3}L_x^4}.
\end{align*} 
Now
\begin{align*}
\indentalign \left\| \left[ \frac{\hat{v}}{\la \sigma \ra^{1/4}} \right]\spcheck \right\|_{L_t^4L_x^2} = \left\| \int_\xi e^{ix\xi} \int_\tau e^{it\tau} \frac{\hat{v}(\xi,\tau)}{\la \tau+\xi\ra^{1/4}} \, d\tau\, d\xi \right\|_{L_t^4L_x^2} \\
&= \left\| \int_\tau e^{it\tau} \frac{ \hat{v}(\xi,\tau)}{\la \tau+\xi \ra^{1/4}} \, d\tau \right\|_{L_t^4L_\xi^2} \\
&\leq \left\| \int_\tau e^{it\tau} \frac{ \hat{v}(\xi,\tau)}{\la \tau+\xi \ra^{1/4}} \, d\tau \right\|_{L_\xi^2L_t^4} \text{ by Minkowskii}\\
&\leq \left\| \int_\tau e^{it\tau} \frac{ \hat{v}(\xi,\tau-\xi)}{\la \tau \ra^{1/4}} \, d\tau \right\|_{L_\xi^2L_t^4}  \\
&\leq \left\| \int_\tau e^{it\tau} \hat{v}(\xi,\tau-\xi) \, d\tau \right\|_{L_\xi^2L_t^2} \text{ by Sobolev}\\
&= \|\hat{v} \|_{L_\xi^2 L_\tau^2} = \|v\|_{L^2}.
\end{align*}
The estimate on 
$$  \left\| \left[ \frac{\hat{v}_2}{\la \sigma_2 \ra^{1/4+}} \right]\spcheck \right\|_{L_t^{8/3}L_x^4}$$
is obtained by interpolating halfway between $\| [\la \sigma_1 \ra^{-d}\hat{v}_1]\spcheck \|_{L_t^4L_x^\infty} \leq \|v_1\|_{L_t^2L_x^2}$ for $d>\frac{1}{2}$ (Strichartz) and $\|[\hat{v}_1]\spcheck \|_{L_t^2L_x^2} = \|v_1\|_{L_t^2L_x^2}$.  This leaves us to estimate
\begin{equation}
\label{E:601}
S' = \int_\ast \frac{\hat{v} \hat{v}_1 \hat{v}_2 | \xi |^{1/2}}{ \la \sigma \ra^b \la \sigma_1 \ra^{c_1} \la \sigma_2 \ra^{b_1}}
\end{equation}
with $\hat{v}=\hat{v}(\xi,\tau)$,  $\hat{v}_1=\hat{v}_1(\xi_1,\tau_1)$, $\hat{v}_2=\hat{v}_2(\xi_2,\tau_2)$, $\sigma_1 = \tau_1 + \xi_1^2$, $\sigma_2=\tau_2+\xi_2^2$, $\sigma = \tau+\xi$, and $*$ indicates the restriction $\xi_1=\xi+\xi_2$, $\tau_1=\tau+\tau_2$. \\

\noindent \underline{Region $\sigma$ dominant, $|\sigma|\geq \max(|\sigma_1|, |\sigma_2|)$.}

\begin{align}
\notag \indentalign S' \leq \left( \int_{\xi,\sigma} |\hat{v}|^2 \right)^{1/2} \left( \int_{\xi,\sigma} \la \sigma \ra^{-2b} \left| \int_{\xi_2,\sigma_2} \frac{ |\hat{v}_1 \hat{v}_2| |\xi|^{1/2}}{ \la \sigma_1 \ra^{c_1} \la \sigma_2 \ra^{b_1}}  \right|^2 \right)^{1/2} \\
\notag &\leq \left( \int_{\xi,\sigma} |\hat{v}|^2 \right)^{1/2} \left( \int_{\xi,\sigma} \la \sigma \ra^{-2b} \left( \int_{\xi_2,\sigma_2} |\hat{v}_1 \hat{v}_2|^2 \right) \left( \int_{\xi_2,\sigma_2} \frac{|\xi|}{\la \sigma \ra^{2c_1} \la \sigma_2 \ra^{2b_1}} \right)  \right)^{1/2} \\
\label{E:A1} &\leq \left(\sup_{\sigma,\xi} \la \sigma \ra^{-2b} \int_{\sigma_2}\int_{\xi_2} \frac{|\xi|}{\la \sigma_1 \ra^{2c_1} \la \sigma_2 \ra^{2b_1}} \, d\xi_2 d\sigma_2 \right)^{1/2} \|v\|_2 \|v_1\|_2 \|v_2\|_2. 
\end{align}
The inner integral is taken over fixed $\sigma, \xi,\sigma_2$.  Since $\xi_1=\xi+\xi_2$, we have $d\xi_1=d\xi_2$, and since $\sigma_1 - \sigma - \sigma_2 = (\xi_1-\frac{1}{2})^2 - ( \xi_2-\frac{1}{2})^2$, we have
$$d\sigma_1 = [+2(\xi_1-\tfrac{1}{2}) - 2(\xi_2-\tfrac{1}{2})] d\xi_1
= 2\xi d\xi_2 . $$
Thus, the quantity in parentheses in \eqref{E:A1} is bounded by
\begin{align*}
\indentalign \sup_\sigma \la \sigma \ra^{-2b} \int_{\sigma_2} \int_{\sigma_1} \frac{d\sigma_1 d\sigma_2}{\la \sigma_1 \ra^{2c_1} \la \sigma_2 \ra^{2b_1}} \\
&\leq \la \sigma \ra^{-2b} \la \sigma \ra^{[1-2c_1]_+} \la \sigma \ra^{[1-2b_1]_+} \\
&\leq \la \sigma \ra^{-2b+[1-2c_1]_+ + [1-2b_1]_+}
\end{align*}
since $\sigma$ is dominant.  If $b,c_1,b_1<\frac{1}{2}$, then the exponent is $2-2b-2c_1-2b_1$, and it suffices to have $b+c_1+b_1\geq 1$.\\

\noindent \underline{Region $\sigma_1$ dominant, $|\sigma_1|\geq \max(|\sigma|,|\sigma_2|)$.}

  By the Cauchy-Schwarz method, it suffices to show
\begin{equation}
\label{E:A2}
\sup_{\xi_1,\sigma_1} \la \sigma_1 \ra^{-2c_1} \int_{\sigma_2} \int_{\xi_2} |\xi| \la \sigma \ra^{-2b} \la \sigma_2 \ra^{-2b_1} \, d\xi_2 d\sigma_2
\end{equation}
is finite.\\
\textit{Subregion} $|\xi_1-\frac{1}{2}|\leq 2|\xi_2-\frac{1}{2}|$.  Then $|\xi|\leq 3|\xi_2-\frac{1}{2}|$.  The inner integral over $\xi_2$ is taken with $\sigma_1$, $\xi_1$, $\sigma_2$ fixed.  Since $\sigma_1-\sigma-\sigma_2=(\xi_1-\frac{1}{2})^2-(\xi_2-\frac{1}{2})^2$, we have $d\sigma = 2(\xi_2-\frac{1}{2})d\xi_2$ and thus \eqref{E:A2} is bounded by
\begin{align*}
\indentalign  \sup_{\sigma_1} \la \sigma_1 \ra^{-2c_1} \int_{\sigma_2=0}^{|\sigma_1|} \la \sigma_2 \ra^{-2b_1} d\sigma_2 \int_{\sigma=0}^{|\sigma_1|} \la \sigma \ra^{-2b} \, d\sigma\\
&\leq \la \sigma_1 \ra^{-2c_1+ [1-2b_1]_+ + [1-2b]_+}.
\end{align*}
If $b,b_1<\frac{1}{2}$, then the exponent here is $2-2b_1-2b-2c_1$, and thus we need $b_1+b+c_1\geq 1$.

\noindent\textit{Subregion} $|\xi_1-\frac{1}{2}|\geq 2|\xi_2-\frac{1}{2}|$.  Since $\xi=\xi_1-\xi_2$, we have $|\xi|\leq \frac{3}{2}|\xi_1-\frac{1}{2}|$.  Also $\frac{3}{4}(\xi_1-\frac{1}{2})^2 \leq (\xi_1-\frac{1}{2})^2-(\xi_2-\frac{1}{2})^2 = \sigma_1-\sigma_2-\sigma \leq 3|\sigma_1|$, and thus $(\xi_1-\frac{1}{2})^2\leq 4|\sigma_1|$.  Thus \eqref{E:A2} is bounded by
$$\sup_{\xi_1,\sigma_1} \la \xi_1 \ra^{1-4c_1}
\int_{\sigma_2}\int_{\xi_2} \la \sigma \ra^{-2b} \la \sigma_2
\ra^{-2b_1} \, d\xi_2 \, d\sigma_2 .$$
We change variables $y=(\xi_2-\frac{1}{2})^2$ to obtain
\begin{equation}
\label{E:A3}
\sup_{\xi_1,\sigma_1} \la \xi_1 \ra^{1-4c_1} \int_{y=-\la \xi_1
  \ra^2}^{y=\la \xi_1 \ra^2} |y|^{-1/2} \int_{\sigma_2} \la \sigma
\ra^{-2b} \la \sigma_2 \ra^{-2b_1} \, d\sigma_2\, dy .
\end{equation}
The inner integral over $\sigma_2$ is taken with fixed $y=(\xi_2-\frac{1}{2})^2$, $\xi_1$, $\sigma_1$, and thus $\sigma_2+\sigma = \sigma_1-(\xi_1-\frac{1}{2})^2-(\xi_2-\frac{1}{2})^2=\sigma_1-(\xi_1-\frac{1}{2})^2-y$ is fixed.  By Lemma \ref{L:GTV4.2},
$$\int_{\sigma_2} \la \sigma \ra^{-2b} \la \sigma_2 \ra^{-2b_1} \, d\sigma_2 \leq c\la \sigma_1-(\xi_1-\tfrac{1}{2})^2-y \ra^{-\alpha} \qquad \text{if }b_1+b>\tfrac{1}{2}$$
where $\alpha = 2b-[1-2b_1]_+$ if $b_1\geq b$ or $=2b_1-[1-2b]_+$ if $b\geq b_1$.  By Lemma \ref{L:GTV4.1} with $f(y)=\chi_{-\la \xi_1 \ra^2\leq y\leq \la \xi_1 \ra^2}(y)|y|^{-1/2}$ and $g(y)=\la y \ra^{-\alpha}$,
\begin{align*}
\indentalign \sup_{\sigma_1} \int_{y=-\la \xi_1 \ra^2}^{+\la\xi_1 \ra^2} |y|^{-1/2} \la \sigma_1-(\xi_1-\tfrac{1}{2})^2-y \ra^{-\alpha} \, dy \\
&\leq  \int_{y=-\la \xi_1 \ra^2}^{+\la \xi_1\ra^2} |y|^{-1/2} \la y \ra^{-\alpha} \, dy \\
&\leq \la \xi_1 \ra^{[1-2\alpha]_+}
\end{align*}
and hence \eqref{E:A3} is controlled by
$$\sup_{\xi_1} \la \xi_1 \ra^{1-4c_1+[1-2\alpha]_+} .$$
We now consider the exponent. Suppose $b,b_1<\frac{1}{2}$ but $b_1+b>\frac{1}{2}$.  Then $\alpha =-1+2b+2b_1$.\\
\textit{Case 1}.  $\alpha>\frac{1}{2} \Longleftrightarrow b+b_1>\frac{3}{4}$.  Then we need $c_1\geq \frac{1}{4}$.\\
\textit{Case 2}.  $\alpha = \frac{1}{2} \Longleftrightarrow b+b_1=\frac{3}{4}$.  Then we need $c_1>\frac{1}{4}$.\\
\textit{Case 3}.  $\alpha<\frac{1}{2} \Longleftrightarrow b+b_1<\frac{3}{4}$.  Then the exponent is $4-4c_1-4b-4b_1$, and we need $b+b_1+c_1\geq 1$.\\

\noindent \underline{Region $\sigma_2$ dominant, $|\sigma_2|\geq \max(|\sigma|,|\sigma_1|)$.}
This is analogous to the $\sigma_1$ dominant case, but we carry it out anyway.  By the Cauchy-Schwarz method, need to show
\begin{equation}
\label{E:A4}
\sup_{\xi_2, \sigma_2} \la \sigma_2 \ra^{-2b_1}
\iint_{\sigma_1,\xi_1} |\xi| \la \sigma \ra^{-2b} \la \sigma_1
\ra^{-2c_1} \, d\sigma_1 \, d\xi_1 
\end{equation}
is finite.\\
\textit{Subregion $|\xi_2-\frac{1}{2}|\leq 2|\xi_1-\frac{1}{2}|$}.  Then $|\xi| \leq 3|\xi_1-\frac{1}{2}|$.  The inner integral over $\xi_1$ is taken with $\sigma_1$, $\xi_2$, $\sigma_2$ fixed.  Since $\sigma_1-\sigma-\sigma_2 = (\xi_1-\frac{1}{2})^2-(\xi_2-\frac{1}{2})^2$,  we have $-d\sigma = 2(\xi_1-\frac{1}{2})d\xi_1$.  Thus \eqref{E:A4} is bounded by
\begin{align*}
\indentalign \sup_{\sigma_2} \la \sigma_2 \ra^{-2b_1} \int_{\sigma_1} \la \sigma_1 \ra^{-2c_1} d\sigma_1 \int_\sigma \la \sigma \ra^{-2b} \, d\sigma\\
&\leq \sup_{\sigma_2} \la \sigma_2 \ra^{-2b_1+[1-2c_1]_+ +[1-2b]_+} .
\end{align*}
If $c_1,b<\frac{1}{2}$, then need $b_1+c_1+b\geq 1$.

\noindent \textit{Subregion $|\xi_2-\frac{1}{2}|\geq 2|\xi_1-\frac{1}{2}|$}.  Then $|\xi| \leq \frac{3}{2}|\xi_2-\frac{1}{2}|$.  Also, $\frac{3}{4}(\xi_2-\frac{1}{2})^2 \leq (\xi_2-\frac{1}{2})^2-(\xi_1-\frac{1}{2})^2 =-\sigma_1+\sigma+\sigma_2\leq 3|\sigma_2|$ and hence $(\xi_2-\frac{1}{2})^2 \leq 4|\sigma_2|$.  We change variable $y=(\xi_1-\frac{1}{2})^2$ to obtain that \eqref{E:A4} is bounded by
\begin{equation}
\label{E:A5}
\sup_{\sigma_2,\xi_2} \la \xi_2 \ra^{1-4b_1} \int_{y=-\la \xi_2
  \ra^2}^{\la \xi_2 \ra^2} |y|^{-1/2} \int_{\sigma_1} \la \sigma
\ra^{-2b}\la \sigma_1 \ra^{-2c_1} \, d\sigma_1 \; dy .
\end{equation}
Since $-\sigma+\sigma_1 = \sigma_2 + (\xi_1-\frac{1}{2})^2-(\xi_2-\frac{1}{2})^2=\sigma_2-(\xi_2-\frac{1}{2})^2+y$ is fixed, by Lemma \ref{L:GTV4.2},
$$\int_{\sigma_1} \la \sigma \ra^{-2b}\la \sigma_1 \ra^{-2c_1} d\sigma_1 \leq \la \sigma_2 - (\xi_2-\tfrac{1}{2})^2+y \ra^{-\alpha} \qquad \text{if }b+c_1>\tfrac{1}{2}$$
with $\alpha=2b-[1-2c_1]_+$ if $c_1\geq b$ and $\alpha=2c-[1-2b]_+$ if $b\geq c_1$.  By Lemma \ref{L:GTV4.1} with $f(y)=\chi_{-\la \xi_2 \ra^2 \leq y \leq \la \xi_2 \ra^2}(y)|y|^{-1/2}$ and $g(y)=\la y \ra^{-\alpha}$,
\begin{align*}
\indentalign \sup_{\sigma_2} \int_{y=-\la \xi_2\ra^2}^{\la \xi_2 \ra^2}
|y|^{-1/2} \la \sigma_2-(\xi_2-\tfrac{1}{2})^2+y\ra^{-\alpha} dy \\
&\leq \int_{y=-\la \xi_2 \ra^2}^{\la \xi_2 \ra^2} |y|^{-1/2}\la y
\ra^{-\alpha} \, dy \\
&\leq \la \xi_2 \ra^{[1-2\alpha]_+}. 
\end{align*}
Hence \eqref{E:A5} is bounded by
$$\sup_{\xi_2} \la \xi_2 \ra^{1-4b_1+[1-2\alpha]_+} .$$
We now consider the exponent. If $b,c_1<\frac{1}{2}$ and $b+c_1>\frac{1}{2}$, then $\alpha=-1+2b+2c_1$.\\
Case 1.  $\alpha<\frac{1}{2} \Longleftrightarrow b+c_1<\frac{3}{4}$.  Then the exponent is $4-4b_1-4b-4c_1$ so we need $b+b_1+c_1\geq 1$.\\
Case 2.  $\alpha=\frac{1}{2} \Longleftrightarrow b+c_1=\frac{3}{4}$.  Here, we need $b_1>\frac{1}{4}$.\\
Case 3.  $\alpha>\frac{1}{2} \Longleftrightarrow b+c_1>\frac{3}{4}$.  Here, we need $b_1\geq \frac{1}{4}$.
\end{proof}

\begin{proof}[Proof of Lemma \ref{L:Zakmultest} \eqref{I:Zakwavecoup}]
We show that the proof of Lemma \ref{L:Zakmultest} \eqref{I:Zakwavecoup} is actually identical to that of Lemma \ref{L:Zakmultest} \eqref{I:ZakSchcoup}.  We discuss only the $+$ case.  The estimate is equivalent to showing
$$|W| \leq c \|v\|_2\|v_1\|_2 \|v_2\|_2$$
where
\begin{equation}
\label{E:610}
W = \int_* \frac{ \hat{v} \hat{v}_1 \hat{v}_2 |\xi| \la \xi \ra^{-1/2}}{\la \sigma \ra^c \la \sigma_1 \ra^{b_1} \la \sigma_2 \ra^{b_1}}
\end{equation}
with $\hat{v}=\hat{v}(\xi,\tau)$, $\hat{v}_1=\hat{v}_1(\xi_1,\tau_1)$, $\hat{v}_2(\xi_2,\tau_2)$, $\sigma = \tau+\xi$, $\sigma_1=\tau_1+\xi_1^2$, $\sigma_2 = \tau_2 + \xi_2^2$, and $*$ indicates the restriction $\xi=\xi_1-\xi_2$, $\tau=\tau_1-\tau_2$. Indeed, for $\hat{v}\in L^2$,
\begin{align*}
\indentalign \int_{\xi,\tau} [ \partial_x(u_1\bar{u}_2)]\sphat(\xi,\tau) \la \xi \ra^{-1/2} \la \sigma \ra^{-c}\hat{v}(\xi,\tau) d\xi d\tau \\
&= \int_{\xi,\tau}  |\xi| \left[ \int_{ \substack{\xi=\xi_1-\xi_2 \\
      \tau=\tau_1-\tau_2}} \hat{u}_1(\xi_1,\tau_1)
  \hat{u}_2(\xi_2,\tau_2) \right] \la \xi \ra^{-1/2} \la \sigma
\ra^{-c}\hat{v}(\xi,\tau) d\xi d\tau .
\end{align*}
Set $\hat{v}_1(\xi_1,\tau_1) = \hat{u}_1(\xi_1,\tau_1)\la \sigma_1 \ra^{b_1}$, $\hat{v}_2(\xi_2,\tau_2) = \hat{u}_2(\xi_2,\tau_2)\la \sigma_2 \ra^{b_1}$ to obtain \eqref{E:610}.

We note that \eqref{E:610} is the same as \eqref{E:600} considered in the proof of Lemma \ref{L:Zakmultest} \eqref{I:ZakSchcoup} with $b$, $c_1$ in \eqref{E:600} replaced by $c$, $b_1$, respectively.  Thus, the condition $b+b_1+c_1\geq 1$ in Lemma \ref{L:Zakmultest} \eqref{I:ZakSchcoup} becomes $c+2b_1\geq 1$ in  Lemma \ref{L:Zakmultest} \eqref{I:Zakwavecoup}
\end{proof}

\newcommand{\etalchar}[1]{$^{#1}$}
\providecommand{\bysame}{\leavevmode\hbox to3em{\hrulefill}\thinspace}
\providecommand{\MR}{\relax\ifhmode\unskip\space\fi MR }
\providecommand{\MRhref}[2]{%
  \href{http://www.ams.org/mathscinet-getitem?mr=#1}{#2}
}
\providecommand{\href}[2]{#2}

\end{document}